\magnification=\magstephalf
\input amstex
\documentstyle{amsppt}
\hoffset 1.25truecm
\hsize=12.4 cm
\vsize=19.7 cm
\TagsOnRight

\font\tenscr=rsfs10 
\font\sevenscr=rsfs7 
\font\fivescr=rsfs5 
\skewchar\tenscr='177 \skewchar\sevenscr='177
\skewchar\fivescr='177
\newfam\scrfam \textfont\scrfam=\tenscr \scriptfont\scrfam=\sevenscr
\scriptscriptfont\scrfam=\fivescr

\newsymbol\leqs 1336
\newsymbol\geqs 133E
\def\le{\leqs}
\def\ge{\geqs}

\def\az{\alpha}  \def\bz{\beta}
    
    \def\fz{\varphi}
  
\def\lz{\lambda}

\def\flsh{\flushpar}

\def\qd{\quad}
\def\qqd{\qquad}

\def\lmm{Lemma}
\def\prp{Proposition}
\def\thm{Theorem}
\def\crl{Corollary}

\def\prf{Proof}

\def\le{\leqslant}
\def\ge{\geqslant}

\topmatter
\title{A Comment on the Book ``Continuous-Time Markov Chains''
by W.J. Anderson} \endtitle
\rightheadtext{A Comment on the Book ``Continuous-Time Markov Chains''}
\author{Mu-Fa Chen}\endauthor
\affil{(Beijing Normal University)}\endaffil
\thanks {Research supported in part by NSFC and the State Education Commission of China.}\endthanks
\address {Department of Mathematics, Beijing Normal University, Beijing 100875,
    The People's Republic of China.} \endaddress
\endkeywords
\abstract{The book ``Continuous-Time Markov Chains'' by W. J. Anderson collects
a large part of the development in the past thirty years. It is now a popular
reference for the researchers on this subject or related fields. Unfortunately,
due to a misunderstanding of the approximating methods, several results in the
book are incorrectly stated or proved. Since the results are related to the
present author's work, to whom it may be a duty to correct the mistakes in
order to avoid further confusion. We emphasize the approximating methods because
they are useful in many situations.}
\endabstract
\endtopmatter

\document
\topinsert
\captionwidth{10 truecm}
\flsh Chin. J. Appl. Prob. Stat. 12:1 (1996), 55-59.
\endinsert

Throughout the note, take $E=\{0, 1, 2, \cdots\}$ and suppose that the
$Q$-matrix $Q=(q_{ij}: i, j\in E)$ in the consideration is totally stable:
$q_i:= - q_{ii}<\infty$ for all $i\in E$. Refer to [1], [3] or [6] for further
notation used below. Here, we mention only that the term ``$Q$-process'' is
called ``$q$-function'' in [1].

\subheading 1 {\bf Stochastic Comparability}.

\medskip

Given two $Q$-matrices $Q^{(1)}= \big(q_{ij}^{(1)} \big)$ and $Q^{(2)}= \big(q_{ij}^{(2)} \big)$,
denote by $P^{\min\,(k)}(t)=\!\big(P_{ij}^{\min\,(k)}(t) \big)$, $k=1, 2$ the corresponding
minimal $Q$-processes. The $Q$-processes $P^{\min\,(1)}(t)$ and
$P^{\min\,(2)}(t)$ are said to be {\bf stochastically comparable}
if
$$ \sum_{j\ge k} P_{ij}^{\min\,(1)}(t) \le  \sum_{j\ge k}
P_{mj}^{\min\,(2)}(t)  \quad  \text{for all $i\le m$ and $k\ge 0$}.  \tag1.1$$
From (1.1), we obtain
$$ \sum_{j\ge k} q_{ij}^{(1)}\le  \sum_{j\ge k} q_{mj}^{(2)}
\quad  \text{for all $i\le m$ and } k\in \{0, \cdots, i\} \cup \{m+1, m+2, \cdots\}.
\tag1.2$$
To see this, simply use the backward Kolmogorov equation to deduce that
$$\lim_{t\to 0}  \frac{1}{t} \sum_{j\in A} P_{ij}^{\min\,(k)}(t)= \sum_{j\in A} q_{ij}^{(k)},
 \quad  \quad i\notin A\subset E. $$
However, the proof of $(1.1)\Longrightarrow (1.2)$ given in [1; p.249] works
only for the conservative case.
The Kirstein's original theorem states that {\it if the $Q$-matrices $Q^{(1)}$ and
$Q^{(2)}$ are both regular, then $(1.1)$ and $(1.2)$ are equivalent}. The theorem
was extended in [1; \thm\;7.3.4], where it was claimed that if both
$Q^{(1)}$ and $Q^{(2)}$ are conservative, then (1.2) implies (1.1) {\it without
using the uniqueness assumption for the $Q$-processes} (By the way, we mentioned
that in the last two sentences of part (1) (resp. part (2)) of \thm\;7.3.4 there,
the items $(a)$ and $(b)$ have to be exchanged). Unfortunately, the extension
is incorrect.

\proclaim{\bf Counterexample} Consider the conservative birth-death $Q$-matrices
$$q_{i, i-1}^{(1)} = q_{i, i+1}^{(1)}=q_{i, i-1}^{(2)}=(i+1)^2, \quad
q_{i, i+1}^{(2)}=\az (i+1)^2,  \quad  i\ge 0,  \quad \az>1.$$
Then $(1.2)$ holds but $(1.1)$ does not.  \endproclaim

\demo{\prf} a) Note that in the conservative case, condition (1.2) is reduced
to that for all $i\le m$,
$$ \sum_{j\ge k} q_{ij}^{(1)}\le  \sum_{j\ge k} q_{mj}^{(2)}
\quad  \text{if $k\ge m+1$} \quad \text{and} \quad
\sum_{j=0}^k q_{ij}^{(1)}\ge  \sum_{j=0}^k q_{mj}^{(2)}
\quad  \text{if $k\le i-1$}. \tag1.3$$
Now, it is easy to check that our birth-death $Q$-matrices satisfy (1.3).

b) Next, for a given birth-death $Q$-matrix: $q_{i, i-1}=(i+1)^2$ and
$q_{i, i+1}=\bz (i+1)^2$, $\bz>0$, it is well known that the $Q$-process
is unique iff $\bz \le 1$. Applying this to our example, we see that
the $Q^{(1)}$-process is unique but not the $Q^{(2)}$-process. Thus,
$ \sum_{j\ge 0} P_{ij}^{\min\,(1)}(t) = 1 >  \sum_{j\ge 0}
P_{mj}^{\min\,(2)}(t)$ and so (1.1) fails. \qed\enddemo

To explain what was wrong in the proof of [1; \thm\;7.3.4], we should introduce
some approximating methods.

\subheading 2 {\bf Approximating Methods}.

\medskip

In the study of Markov chains, there are several different approximating
methods. Among them, the simplest one is as follows. Take $E_n\subset E$,
$E_n\uparrow E$ and let $\{Q_n\}$ be the truncated $Q$-matrices of $Q$. That
is,
$$q_{ij}^{(n)}=
\cases
q_{ij},  \quad &\text{if } i, j \in E_n\\
0, &\text{otherwise} \endcases    \tag2.1 $$
Then we have
$$P_{ij}^{\min\, (n)}(t) \uparrow P_{ij}^{\min}(t)
 \quad \text{as } n\uparrow\infty\; \; \text{for all $i, j\in E_n$ and $t\ge 0$}.
\tag2.2$$
Usually, one chooses $\{E_n\}$ so that each $Q_n$ is a bounded $Q$-matrix.
Note that
the $Q$-matrices $Q_n$ often become non-conservative even though the original
$Q$-matrix is usually assumed to be conservative in practice.

It is known that for bounded $Q$-matrices, (1.1) and (1.2) are equivalent.
Hence, the main step of the proof of [1; \thm\;7.3.4] is again using a sequence of $Q$-processes with
bounded $Q$-matrices to approximate the minimal $Q$-process. However, in order
to keep (1.2), the truncated $Q$-matrices given by (2.1) is not suitable. One
adopts a different choice:
$$q_{ij}^{(n)}=
\cases
q_{ij},  \quad  &\text{if $i, j \le n-1$}\\
 \sum_{k\ge n} q_{ik},  &\text{if $i\le n-1$ and $j=n$}\\
0, &\text{otherwise.}\endcases
\tag2.3$$

Two different approximation methods were introduced in the study of
reaction-diffusion processes [2, 3, 4]. The first one is a truncating from $E$
to $E_n:=\{i\in E: q_i\le n\}$,
$$q_{ij}^{(n)}=I_{E_n}(i) q_{ij}, \qqd i,\; j\in E,  \tag2.4$$
where $I_A$ is the indicator of set $A$. The second one is stopped at the
$n$-th row,
$$q_{ij}^{(n)}=
\cases
q_{ij},\qqd &\text{if} \qd i\le n\\
q_{nj},\qqd &\text{if} \qd i> n.
\endcases\tag 2.5$$
In both cases, the resulting $Q$-matrices $Q_n=\big(q_{ij}^{(n)}\big)$ are conservative if
so is the original $Q=(q_{ij})$. The latter one was used to keep the Lipschitz
property of the corresponding semi-group [2]. The former one is especially
powerful to deal with the uniqueness of $Q$-processes. It was originally used
in [3] to prove [1; \crl\;2.2.15, \crl\;2.2.16 and \thm\;7.5.5].
It should be pointed out that the results of [1; \S 2.2] starting from \prp\;2.2.12
are mainly taken (but without mentioned) from [3] restricted to the particular
case of Markov chains.

From my knowledge, a new point was made by the author
in [1; \prp\;2.2.14], where (2.1), (2.3) and (2.4) are unified into a general
form. Take
$E_n\uparrow E $ and let $Q_n=\big(q_{ij}^{(n)}\big)$ be a $Q$-matrix such that
$$q_{ij}^{(n)}=
\cases
q_{ij},\qqd &\text{if}\qd i,\; j\in E_n\\
0,\qqd &\text{if}\qd i\notin E_n.
\endcases\tag 2.6$$
Here $q_{ij}^{(n)}\, (i\in E_n,\; j\notin E_n)$ are left to be freedom\footnote{Addition in Proof (This and the next footnotes
are made on December 24, 2014). Certainly, the additional conditions
that $q_{ij}^{(n)}\ge 0$ for all $i\in E_n, j\ne E_n$ and $\sum_{j}q_{ij}^{(n)} \le 0$ for every $i\in E_n$ are required here.}.
Moreover, the author proved that
$$\align
&P_{ij}^{\min\, (n)}(t)\le P_{ij}^{\min\, (n+1)}(t)
  \quad \quad \text{for all } i,\; j\in E_n  \; \text{and $t\ge 0$}\\
&\text{and}\; \lim_{n\to\infty} P_{ij}^{\min\, (n)}(t) = P_{ij}^{\min}(t)
 \quad  \quad  \text{for all $i, j$ and $t\ge 0$}.
\tag2.7  \endalign$$
Since the monotonicity of $P_{ij}^{\min\, (n)}(t)$ in $n$ depends on $i$ and
$j$, the conclusion (2.7) is weaker than the following statement:
$$ P_{ij}^{\min\, (n)}(t) \uparrow P_{ij}^{\min}(t)  \quad \text{as } \;
n\to \infty, \quad \text{for all }i, j, t\ge 0.$$
Unfortunately, the author neglected the difference of these two statements
and incorrectly wrote the latter one as the conclusion of [1; \prp\;2.2.14].\footnote{Here are some
counterexamples to the proposition.
Let $Q=(q_{ij}: i, j\in E)$ be a conservative matrix on $E=\{0, 1, 2, \ldots\}$.
For each $n\ge 1$, set $E_n=\{0, 1, \ldots, n\}$ and define
$$q_{ij}^{(n)}= I_{E_n}(i) q_{ij}, \qquad i, j \in E$$
which is (2.4) here and is a special case of (2.6) in [1].
This is a bounded and conservative $Q$-matrix on $E$ and hence the corresponding process is unique.
For each $n$, denote  by $P^{(n)}(t)$ the corresponding process. Of course,
$$\sum_{j\in E} P_{ij}^{(n)}(t)=1,\qquad i\in E,\; t\ge 0.$$
Recall that [1; Proposition 2.2.14] says that
$$P_{ij}^{(n)}(t)\uparrow P_{ij}^{\min}(t), \qquad i,j \in E,\; t\ge 0$$
where $P^{\min}(t)$ is the minimal process determined by the original $Q$.
If the proposition were true, then the monotone convergence theorem would imply that
$$1=\sum_{j\in E} P_{ij}^{(n)}(t) \uparrow \sum_{j\in E} P_{ij}^{\min)}(t)\qquad \text{as }n\to\infty.$$
This leads to a contradiction since we do not assume here the $Q$-process to be unique.\newline
\text{\quad} As far as I know, among the different approximating methods listed at the beginning of this
section, the conclusion of [1; Proposition 2.2.14] holds only for the one defined by (2.1).\newline
\text{\quad} Clearly, the proof of [1; Corollary 2.2.15] is incomplete since at the last step, one uses the incorrect [1; Proposition 2.2.14].\newline
\text{\quad} However, [1; Corollary 2.2.16] is correct. Honestly, I would to say that
this is one of my favourite
contribution to the theory of Markov Chains. It appeared first in
my Chinese book [4] in 1986 and then in the paper [3; Theorem 16)]
in the same year. This paper is cited in [1], but the originality of
the result is not mentioned.
As we know, the earlier known criterion for the uniqueness is asking to solve
a equation with infinite variables. This is usually not
practicable. Actually, it costs me more than 5 years to find
out such a powerful sufficient condition. A short story
about this and others is included in Chapter 9 of my book
``{\it Eigenvalues, Inequalities, and Ergodic Theory}''.
Springer 2005.
In which, a number of papers are included for the applications
of this result. A new application of the result to genetic
study is given in ``Transformation Markov jump processes'' by Z.M. Ma et al.
There the authors apply the result to the continuous state space.
Actually, my original result is stated for general state
space, not necessarily discrete one.
Except the 1986's publications mentioned above, the
same result was also published several times later:
in the paper ``{\it On three classical problems for Markov chains with continuous
      time parameters}, J. Appl. Prob. 28(1991)2, 305-320;
in the enlarged version [6; pages 81 and 84] (2nd ed. 2004) of my Chinese book [4];
in a textbook ``{\it Introduction to Stochastic Processes} (in Chinese, 2006.
Higher Education Press, Beijing)'' by me and Yong-Hua Mao.
Actually, this result was taught in a course on Stochastic Processes for
undergraduate/graduate students at my university every year since 1989.}
This leads to the incorrect proof of [1; \thm\;7.3.4]. More precisely,
in the latter case, we have
$$\lim_{n\to\infty}\sum_{j\ge k} P_{ij}^{\min\, (n)}(t)=
\sum_{j\ge k} P_{ij}^{\min}(t).$$
But under (2.7), we may only have
$$\varliminf_{n\to\infty}\sum_{j\ge k} P_{ij}^{\min\, (n)}(t)\ge
\sum_{j\ge k} P_{ij}^{\min}(t).$$
The inequality can be appeared as shown by the above counterexample.

Note that (2.5) is not contained in (2.6). In general if
$$\lim_{n\to\infty} q_{ij}^{(n)}=q_{ij}, \quad  \quad i, j\in E,\tag2.8$$
then by the backward Kolmogorov equation and the Fatou's lemma,
we have
$$\varliminf_{n\to\infty} P_{ij}^{\min\,(n)}(t)\ge P_{ij}^{\min}(t),
\quad  \quad i, j\in E, \; \; t\ge 0\tag2.9$$
(cf. [6; \lmm\;5.14] for example). For such a general setup, the conclusion
(2.9) is certainly far away from (2.7),\,we even do not know that the limit\,$\lim_{n\to\infty}\!P_{ij}^{\min (n)}\!(t)$
exists or not and can not say that
$\varliminf_{n\to\infty} P_{ij}^{\min\,(n)}(t) $
provides us a $Q$-process. However, whenever $Q=(q_{ij})$ being regular, we
do have
$$\lim_{n\to\infty} P_{ij}^{\min\,(n)}(t)= P_{ij}^{\min}(t),
\quad  \quad i, j\in E, \; \; t\ge 0.$$

\subheading 3 {\bf Extended Kirstein's Theorem}.

\medskip

The correct extension of the Kirstein's Theorem may be stated
as follows.

\proclaim {\bf \thm} The condition $(1.1)$ always implies $(1.2)$. Conversely,
$(1.2)$ implies $(1.1)$ provided either $Q^{(1)}$ and $Q^{(2)}$ are bounded
or $Q^{(2)}$ being regular.
\endproclaim

Because of (2.7) or (2.9), the $Q$-matrix $Q^{(1)}$ can be arbitrary. The
proof of the conclusion is the
same as the proof of [5; \thm\;7] or [6; \thm\;5.31] and hence is omitted
here.
By the way, we mention that for [5; \thm\;7] or [6; \thm\;5.31],
the regularity assumption for the first $q$-pair can be removed.

The next result shows that in general it is impossible to remove the uniqueness
assumption. Recall that a $Q$-process $P(t)=(P_{ij}(t))$
is said to be {\bf monotone} if (1.1) holds for the same process $P(t)$.

\proclaim{\bf \prp} Given a single birth $Q$-matrix $($i.e., $q_{i, i+1}>0$ and
$q_{ij}=0 $ for all $j>i+1$, $i\ge 0)$, the corresponding minimal
$Q$-process is monotone iff the $Q$-process is unique and $(1.2)$ holds.
\endproclaim

\demo{\prf} The sufficiency as well as the necessity of (1.2) follow from the
above theorem. To prove the necessity of the uniqueness, suppose that (1.1)
holds. Then, $\sum_{j} P_{ij}^{\min}(t)$ is increasing in $i$ and so is its
Laplace transform $\sum_{j} P_{ij}^{\min}(\lz)\,(\lz>0)$. Hence
$$z_i(\lz):=1- \lz \text{$\tsize\sum_{j}$} P_{ij}^{\min}(\lz),\qqd \lz>0$$
is decreasing in $i$. On the other hand, it is known that $z_i(\lz)$ is increasing in
$i$ and it is indeed strictly increasing except $z_i(\lz)\equiv 0$
[6; Proof of \thm\;3.16]. Therefore, $z_i(\lz)\equiv z_0(\lz)=0$ and hence the
process is unique.\qed\enddemo

The above proposition shows that the example given in [1; p.251] and furthermore
[1; \crl\;7.4.3] are also incorrect.

\subheading 4 {\bf Application to the Uniqueness of $Q$-processes}.

\medskip

To illustrate an application of the comparison technique, we present the
following result.

\proclaim{\bf \prp} Let $Q^{(1)}$ and $Q^{(2)}$ be conservative $Q$-matrices
satisfying $(1.2)$. If $Q^{(2)}$ is regular, then so is $Q^{(1)}$. In particular,
if there exists a non-negative function $\fz$ such that
$\fz_i\uparrow \infty$ as $i\uparrow\infty$ and moreover
$$ \sum_{j} q_{ij}^{(2)} (\fz_j-\fz_i)\le c (1+\fz_i),  \quad  \quad i\in E
\tag4.1$$
for some constant $c\ge 0$. Then both $Q^{(1)}$ and $Q^{(2)}$ are regular.
\endproclaim

\demo{\prf} a) By using (2.3) and restricting to the finite space $E_n=\{0, 1, \cdots, n\}$,
we obtain
$$P_{in}^{(n,1)}(t)\le P_{in}^{(n,2)}(t), \qqd i\le n-1, \qd t\ge 0,
$$
where $\big(P_{ij}^{(n,k)}(t)\big)$ is the $Q$-process corresponding the
truncating $Q$-matrix $Q_n^{(k)}$, $k=1,\;2$.
This means that
$$\sum_{j\ge n}P_{ij}^{(n,1)}(t)\le \sum_{j\ge n}P_{ij}^{(n,2)}(t),
\qqd i\le n-1, \qd t\ge 0. \tag4.2
$$
Next, denote by $\big(\bar P_{ij}^{(n,k)}(t)\big)$ $(k=1,\;2)$ the $Q$-process obtained by
using (2.4). Then, by the localization theorem [6; \thm\;2.13], we get
$$\sum_{j\le n-1} P_{ij}^{(n,k)}(t)
=\sum_{j\le n-1} \bar P_{ij}^{(n,k)}(t), \qqd i\le n-1, \qd k=1,\;2.$$
Equivalently,
$$\sum_{j\ge n} P_{ij}^{(n,k)}(t)
=\sum_{j\ge n} \bar P_{ij}^{(n,k)}(t), \qqd i\le n-1, \qd k=1,\;2.$$
Combining this with (4.2), we obtain
$$\sum_{j\ge n}\bar P_{ij}^{(n,1)}(t)\le \sum_{j\ge n}\bar P_{ij}^{(n,2)}(t),
\qqd i\le n-1, \qd t\ge 0.
$$
The first assertion now follows from [1; \crl\;2.2.15].

b) Let (4.1) holds. Applying [1; \crl\;2.2.16] to $E_n=\{i: \fz_i\le n\}$ and
$x_i=1+\fz_i$, it follows that the $Q^{(2)}$-process is unique.
Hence the second assertion follows from the first one.\qed\enddemo

We remark that the condition (4.1) is necessary for the uniqueness of single
birth processes [3; Remark 23]. Refer to [4; Remark 3.5.1] or [6; Remark 3.20]
for a proof.

Finally, the proof of [1; \thm\;7.5.5] is also incorrect since the resulting
bounded $Q$-matrices do not possess the coupling relation. The original proof
uses (2.4). Refer to [3] or [4; Chapter 5] for details.

\medskip

\flsh{\bf Acknowledgement}. The author is greatly indebted to Prof. Yong-Long
Dai, Prof. Feng-Yu Wang and Mr. Yu-Hui Zhang for their comments on the first
version of the paper. The proof a) of the last proposition is due to Y. H.
Zhang, which improves the original conclusion.

\subheading 5 {\bf References}.

\medskip

\widestnumber\key{20}
\ref\key {}1
\by Anderson, W. J. (1991)
\book Continuous-Time Markov Chains
\publ Springer Series in Statistics  \endref

\ref\key 2
\by Chen, M. F. (1985)
\paper Infinite-dimensional reaction-diffusion processes
\jour Acta Math. Sin. New Ser.{\bf 1}:3, 261--273\endref

\ref\key 3
\by Chen, M. F.(1986)
\paper Couplings of jump processes
\jour Acta Math. Sin. New Ser.{\bf 2}:2, 123--136\endref

\ref\key 4
\by Chen, M. F.(1986)
\book Jump Processes and Particle Systems  \rm{(In Chinese)}
\publ Beijing Normal Univ. Press
\publaddr Beijing\endref

\ref\key 5
\by Chen, M. F.(1991)
\paper  On coupling of jump processes
\jour Chin. Ann. Math. {\bf 12(B)}:4, 385--399\endref

\ref\key 6
\by Chen, M. F. (1992)
\book From Markov Chains to Non-Equilibrium Particle Systems
\publ  World Scientific \endref

\ref\key 7
\by Kirstein, B. M.(1976)
\paper Monotonicity and comparability of time-homogeneous Markov processes
with discrete state space
\jour Math. Operationsforsch. Statist. {\bf 7}, 151--168   \endref
\enddocument